\documentclass[12pt]{amsart}
\usepackage{yhmath}

\usepackage{amssymb,amsmath,amsfonts,amscd,amsthm,enumerate}
 \usepackage{graphicx,color}
\graphicspath{{picturesfolder/}}
\usepackage{wrapfig}
\usepackage{multirow} 
\usepackage{textcomp}

\usepackage{todonotes}

\usepackage{epsfig,subfigure}
\usepackage{todonotes}

\usepackage{lineno,xcolor}
\usepackage[all]{xy}
\usepackage{color}

\newcommand{\R}{{\Bbb R}}
\newcommand{\C}{{\Bbb C}}

\def\CP{{\Bbb C  \Bbb P}}
\def\RP{{\Bbb R  \Bbb P}}
\def\bP{{\Bbb P}}

\newcommand{\s}{\sigma}

\newcommand{\dd}{\displaystyle}
\newcommand{\del}{\partial}

\newtheorem{theorem}{Theorem}

\title{No Hyperbolic Pants for the 4-body Problem}

\author[Connor Jackman]{Connor Jackman}
\address[Jackman]{Mathematics Department, University of California, 
4111 McHenry
Santa Cruz, CA 95064, USA}
\email{cfjackma@ucsc.edu}

\author[Richard Montgomery]{Richard Montgomery}
\address[Montgomery]{Mathematics Department, University of California, 
4111 McHenry
Santa Cruz, CA 95064, USA}
\email{rmont@ucsc.edu}
\begin{document}

\date{\today}
  
\begin{abstract}  The $N$-body problem with a $1/r^2$ potential has, in addition to translation and rotational
symmetry, an effective scale symmetry which allows its zero energy flow  to be reduced to a geodesic flow on complex projective
$N-2$-space, minus a hyperplane arrangement.  When $N=3$  we get a geodesic flow on  the two-sphere minus three points.
If,  in addition we assume that  the three masses are equal, then it was proved in \cite{pants} that the corresponding metric is hyperbolic: its     Gaussian curvature is   negative except at two points.
Does the negative curvature property persist for $N=4$, that is,  in  the equal mass $1/r^2$ 4-body problem?
Here we prove `no' by computing that  the  corresponding Riemannian metric in this $N=4$ case has
  positive sectional curvature
at some   two-planes. This `no' answer    dashes     hopes
of  naively extending  hyperbolicity  from $N=3$ to $N > 3$.
\end{abstract}

\maketitle

%\renewcommand{\thesubsection}{\arabic{subsection}}

%\numberwithin{equation}
\section{Introduction}

In \cite{pants} it was shown that the reduced Jacobi-Maupertuis metric for a certain three-body problem 
had negative Gaussian curvature (except at two points where it is zero).   This  hyperbolicity  led to deep dynamical consequences.
Does  hyperbolicity, i.e. curvature negativity,  persist for the analogous  $N$-body problem, $N >3$? 
No.  We show that the analogous  reduced  4-body problem with its metric has two-planes at which the sectional curvature is positive.

The $N$-body problem in question has  equal masses and the  inverse {\it cube} law attractive force between bodies.

\section{Set-up}  Identify the complex number line $\C$   with the Euclidean plane $\R^2$.Then the planar $N$-body problem has configuration space $\C^N \setminus \Delta$.
 Here $\Delta$ is the ``fat diagonal'' consisting
of all collisions:  $\Delta = \{q = (q_1, q_2, \ldots, q_N) \in \C^N:  q_i = q_j \text{ for some pair } i \ne j \} $.
The quotient  of  $\C^N \setminus \Delta$ by translations and  rotations is the ``reduced $N$-body configuration space'': 
$$C_N = Y_N \times \R^+,  Y_N = \CP^{N-2} \setminus \bP \Delta$$
where $\C \bP^{N-2}$ is the projectivization of the center of mass subspace $\C^{N-1}= \{q \in \C^N: \Sigma m_i q_i = 0 \}$
and $\bP \Delta \subset \CP^{N-2}$ is the projectivization of $\Delta \cap \C^{N-1}$.  The $\R^+$ factor records the overall scale of the planar $N$-gon and 
is coordinatized by $\sqrt I$ with $I = \Sigma m_i |q_i|^2$ being the total moment of inertia about the center of mass.  $Y_N$ is the moduli space of
oriented similarity classes of non-collision $N$-gons and will be called ``shape space.''

The following considerations reduce the zero angular momentum, zero energy $N$-body problem to a geodesic flow on 
shape space $Y_N$, {\it provided} the potential $V$  is homogeneous of degree $-2$.  If $V$ is homogeneous of degree $-\alpha$
then the virial identity, also known as the Lagrange-Jacobi identity, asserts that along solutions of energy $H$ we have 
$\ddot I = 4H - (4-2 \alpha)V$ which implies that the only case in which we can generally guarantee that $\ddot I = 0$ is when
$\alpha =2$ and $H=0$. If in addition $\dot I = 0$ then solutions lie on constant levels of $I$.  
Now  we recall  the Jacobi-Maupertuis [JM] reformulation of mechanics  which asserts that 
the solutions to Newton's equations at energy $H$ are, after a time reparameterization,
precisely the geodesic equations for the {\it Jacobi-Maupertuis metric}
 $$ds^2 _{JM} = 2(H-V) ds^2$$
 on the {\it Hill region} $\{ H - V \ge 0 \} \subset \C^N \setminus \Delta$ with $ds^2$ the mass metric. 
 We are interested in the case $H=0$, $-V > 0$ with  $V$ homogeneous of degree $-2$,
 in which case the Hill region is all of $\C^N \setminus \Delta$ and
 $$ds^2_{JM} = U ds^2  , \qquad U = -V$$
 
 The case of prime interest to us  is 
 \begin{equation}
 \label{potential}
 U = - V = \Sigma m_i m_j / r_{ij} ^{2} 
 \end{equation}
 where the sum is over all distinct pairs $ij$.
 This  $U$, and hence the JM metric,   is invariant under rotations and translations. Quotienting first by translations we take representatives in the totally geodesic center of mass zero subspace $\C^{N-1}$, which reduces the dynamics to geodesics of the metric $ds^2_{JM}|_{\C^{N-1}}$ on $\C^{N-1}$.
 Moreover, $ds^2_{JM}|_{\C^{N-1}}$ is also invariant under scaling since the homogeneities of $U$ and
 the  Euclidean  mass metric $ds^2$ on $\C^{N-1}$ cancel.  Thus the JM metric admits  the group $G=\C^*$ of  rotations and scalings
 as an isometry  group.  Now $Y_N $ is the quotient space: $Y_N=  (\C^{N-1} \setminus \Delta)/G=\CP^{N-2}\backslash \Delta$.  
 (By abuse of notation, we continue to use the symbol $\Delta$ to  denote   the image of the collision locus $\Delta$  under projectivization and intersection.) Insisting that the quotient
 map $\pi: \C^{N-1} \setminus \Delta \to Y_N$ is a Riemannian submersion induces a metric on $Y_N$.
 Recall that this means that we  can define the metric on  $Y_N$   by {\it isometrically}  identifying the tangent space to $Y_N$ at a point  $p$ with the orthogonal
 complement (relative to $ds^2 _{JM}$ or $ds^2$, and at any point lying over $p$ in $\C^{N-1}$) to the $G$-orbit that corresponds to that point.
 These orthogonality conditions are equivalent to the conditions that the linear momentum, angular momentum,
 and `scale momentum' $\dot I$ are all zero.    To summarize, by  using the JM metric and forming the Riemannian quotient,
  the zero-angular momenetum, zero energy
 $1/r^2$ $N$-body problem becomes equivalent to the  problem of finding geodesics for the  metric defined by Riemannian submersion  on $Y_N$. 
 
 {\bf Remark.} The metric quotient procedure just described  realizes the Marsden-Weinstein symplectic reduced space
 of $T^* (\C^N \setminus \Delta)$ by the action of translations, rotations and scalings, $\C\rtimes \C^*$, at momenta values $0$, together with the $N$-body reduced  Hamiltonian flow, 
 but valid only at zero energy.

 {\bf Remark} This  metric  on $Y_N$ can be expressed as $U ds^2_{FS}$ where
 $ds^2 _{FS}$ is the usual Fubini-Study metric on $\CP^{N-2}$.
 
 {\bf Remark} For the standard $1/r^2$ potential of (eq. \ref{potential}) this metric on $Y_N$ is complete, infinite volume.
 
 The collinear $N$-body problem defines a totally geodesic submanifold
 $$\RP ^{N-2} \setminus   \Delta  \subset \CP^{N-2} \setminus  \Delta$$
 We   obtain this   submanifold by  placing the $N$-masses  anywhere along the real axis $\R \subset \C$, arranged
 so their center of mass is zero and so that there are no collisions, and then taking the quotient.
 In other words, $\RP ^{N-2} \setminus   \Delta$ is the  quotient  of $\R^{N-1} \subset \C^{N-1}$ by dilations and real reflections.
 
  \section{Main result}  In case $N=3$, with the potential  (eq. \ref{potential}) above,  $Y_3$ is a pair of pants - a sphere minus three points. The point of \cite{pants} was to show that the metric
 on $Y_3$ just described is hyperbolic  provided $m_1 = m_2 = m_3$.  Specifically, in this equal mass case the
 Gaussian curvature of the metric on the surface $Y_3$ is negative everywhere except at two points (these being the
 ``Lagrange points'' corresponding to equilateral triangles.) What about $Y_4$?
 
\begin{theorem}  Consider the Jacobi-Maupertuis metric on $Y_4$ induced as above for  the case of 4 equal masses under the strong force $1/r^2$ potential  (eq. \ref{potential}). 
  Then  there are   two-planes $\sigma$ tangent to $Y_4$ at which the     Riemannian sectional curvature $\mathcal{K}(\sigma)$ is  positive.

\end{theorem}
 
 {\bf Remark 1.}  The  two-planes $\sigma$ of the theorem pass through    special  points $p \in \RP^2   \subset \CP^2$ which represent certain special collinear configurations. 
See figure \ref{configs}.
 The  two-plane $\sigma$   at  $p$ will be the   orthogonal complement to   $T_p \RP^2$, the normal 2-plane, and is
   realized as   $\sigma = i T_p \RP^2$, using the standard complex structure
 on $\CP^2$. 

{\bf Remark 2.} [Negative curvatures] 
 The $\RP^2$ of  the previous remark  is a totally geodesic surface fixed by an isometric involution.
 There are   other such totally  geodesic surfaces defined as fixed loci of  symmetries, and computer
experiments suggest that these  all    have negative Gaussian curvature everywhere
while their normal 2-planes can have  positive sectional curvature  at some points,   like our special case   $\RP^2$.  Computer experiments also indicate that in the direction of the normal plane there is positive sectional curvature over all collinear configurations of $\RP^2$ and not just the special configurations verified in the theorem. 
An analytic proof of these  claims  beyond our   special case   however looks frightening.

 {\bf Open Question.} A geodesic flow can still be hyperbolic as a  flow, without
 the underlying metric  having all sectional curvatures negative.   Is    geodesic flow on $Y_4$  hyperbolic  as a flow? 
 Is it  even partially hyperbolic?

 \section{Proof of the theorem}

 We take the case $N=4$ in the above considerations.  
 When all the masses are equal to 1 then the mass metric, used to compute the kinetic energy and moment of inertia,
 is the standard Hermitian metric in coordinate $(q_1, q_2, q_3, q_4) \in \C^4$, where the $q_i$ represent the postions
 of the $i$th body. 
 We reduce by translations by going to the  center-of-mass-zero  space 
 which is a 3-dimensional subspace  $\C^3 \subset \C^4$ having Jacobi coordinates as  Hermitian orthonormal  coordinates : 
 
 $$\C^3\stackrel{L}{\to} \C^4 \text{ given by matrix } \begin{bmatrix}\frac12 &\frac{1}{\sqrt 2} & 0\\ \frac12& -\frac{1}{\sqrt 2} & 0\\ -\frac12&0&\frac{1}{\sqrt 2}\\ -\frac12&0&-\frac{1}{\sqrt 2}\end{bmatrix} \text{ in standard bases.}$$

As is well-known, if we start tangent to the center-of-mass-zero subspace  $L(\C^3)$  we stay tangent to it.
Hence we can restrict the dynamics, potential, metric, etc. to the center-of-mass zero subspace. 
We denote the potential restricted to the center of mass zero subspace in Jacobi coordinates as $U_L=U\circ L$
 and still write $ds^2_{JM}=U_Lds^2$ for the restricted JM metric on $\C^3\backslash \Delta$ where 
$ds^2$ is the standard metric on $\C^3$. 

Continuing along the outline above, we now quotient by scaling and rotation isometries, $\C^*$, of $ds^2_{JM}$ to obtain the ``shape space" $Y_4$ and we label the quotient map $\pi: \C^3\backslash \Delta\to Y_4$ which takes a configuration $q$ to it's orbit $\C^*q$. We denote the vertical and horizontal distributions as $\mathcal{V}_p=\ker d_p\pi=\C p$ and $\mathcal{H}_p=\mathcal{V}_p^\perp\stackrel{d\pi}{\cong} T_{\pi(p)}Y_4$.  Requiring $d\pi|_{\mathcal{H}, ds^2_{JM}|_{\mathcal{H}}}$ to be an isometry defines our induced metric on $Y_4$ whose geodesics correspond to $N$-body motions in ``shape space". Under this induced metric on $Y_4$ we denote sectional curvature through the plane $\s\in T_{\pi(p)}Y_4$ by $\mathcal{K}(\s)$. 

Suppressing the notation of evaluating at a representative $p\in \pi(q)$, our main tool in the computation of $\mathcal{K}(\s)$, the $ds^2 _{JM}$ curvature,  is the equation:
 \begin{equation}
 \label{curv}
U_L^3\mathcal{K}(\s)=\frac34((\del_1 U_L)^2+(\del_2 U_L)^2)-\|\nabla U/2\|^2- U_L/2(\del_1^2 U_L+\del_2^2 U_L)+3\frac{U_L^2}{\| p\|^2}(v_1\cdot i v_2)^2
 \end{equation}

Here $\del_a f$ denotes $df(v_a)$ where $f\in C^\infty(\C^3)$ and where  $a =1, 2$ with $v_1, v_2\in \mathcal{H}$ being $ds^2$-orthonormal vectors whose pushforwards $d\pi v_a$ span $\s$. The $'\cdot', \|~\|, \nabla$   refer to
the norm, metric, and Levi-Civita connection for the Euclidean metric
 $ds^2$.  For the derivation of (eq. \ref{curv}) see Appendix A. 

\begin{figure}[ht!]
\centering
\includegraphics[width=90mm] {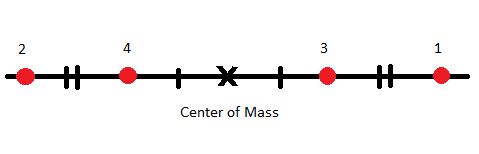}
\caption{: The collinear configurations  $p$ which we consider.}
\label{configs}
\end{figure}

The collinear configurations form a totally geodesic $\RP^2 \subset \CP^2$
which is the image under the projection $\pi$ of the real 2-sphere in $\C^3$  whose points we parameterize as 
 $$p=(\cos\phi\cos\theta, \cos\phi\sin\theta, \sin\phi).$$ We   evaluate (eq. \ref{curv})
 and find positive sectional curvature over the configurations with $\theta=\pi/2$ (see figure \ref{configs}) in the direction of the $iT\RP^2$ plane.
 This plane is  spanned by the pushforwards of
$$v_1=-i\frac{\del p}{\del \phi}=i(\sin\phi\cos\theta, \sin\phi\sin\theta, -\cos\phi)$$
$$v_2=\frac{i}{\cos\phi} \frac{\del p}{\del \theta}=i(-\sin\theta, \cos\theta, 0).$$

TERMS 1: \textit{Over $\RP^2$ in the $iT\RP^2$ direction, the  last term $'\cdot'$ and first  two terms (the first partials)  of (eq. \ref{curv}) vanish:}
$$v_1\cdot iv_2=0,~\del_a U_L=0.$$

That $v_1\cdot iv_2=0$ is clear: $i$ rotates $v_2$ into purely real coordinates. To evaluate the 1st partials, note $Lp$ has purely real coordinates and $\nabla U$ has $k$th component $\sum_{j\neq k}\frac{q_j-q_k}{r_{jk}^4}$, so $\nabla |_{Lp} U$ has purely real coordinates. Now since $Lv_a$ has purely complex coordinates, $$\del_a U_L=\nabla |_{Lp} U\cdot L v_a=0. \qed$$\\

TERMS 2: \textit{With the notation $Lp=(q_1, q_2, q_3, q_4)$, $Lv_a=i(v_a^1, v_a^2, v_a^3, v_a^4)$ and  $\rho_{jk}=\frac{1}{q_j-q_k} ,~\alpha_{jk}=(v_1^j-v_1^k)^2+(v_2^j-v_2^k)^2 \in \R$,  
 the 2nd partials terms of (eq. \ref{curv}) are:}

$$\del_1^2U_L+\del_2^2 U_L=-2\sum_{j>k} \alpha_{jk}\rho_{jk}^4.$$

We write our standard coordinates on $\C^4$ as $q_j=x_j+iy_j$, then since $Lv_a$ is purely imaginary:
$$\del_a^2 U_L=\nabla |_{Lp}(\nabla U\cdot Lv_a)\cdot L v_a=(\nabla |_{Lp} \frac{\del U}{\del y_k}v_a^k)\cdot Lv_a=\frac{\del^2 U}{\del y_j\del y_k}|_{Lp} v_a^kv_a^j.$$
Next we compute $\dd{\frac{\del^2 U}{\del y_j\del y_k}|_{Lp}=2\rho_{jk}^4}$ for $j\neq k$ and $\dd{\frac{\del^2 U}{\del y_k^2}|_{Lp}=-2\sum_{j\neq k} \rho_{jk}^4}$ so now:

$$\del_a^2 U_L=-2\sum_{j\neq k} \rho_{jk}^4((v_a^k)^2-v_a^jv_a^k)=-2\sum_{j>k}\rho_{jk}^4((v_a^k)^2-2v_a^kv_a^j+(v_a^j)^2)=\\ -2\sum_{j>k}\rho_{jk}^4(v_a^k-v_a^j)^2.$$ $\qed$\\

RESULT: \textit{Over the circle $\theta =\pi/2$, $\mathcal{K}(iT\RP^2)$ is positive.}\\

For, substituting terms 1 and 2 into formula (eq. \ref{curv}), we see that: 

$$0<\mathcal{K}\iff 0<U_L^3\mathcal{K}=-\|\nabla U/2\|^2+U_L \sum_{j>k}\alpha_{jk}\rho_{jk}^4\iff $$ 
\begin{equation}
 \label{ineq}
\sum_k(\sum_{j\neq k} \rho_{jk}^3)^2<(\sum_{j>k}\rho_{jk}^2)(\sum_{j>k}\alpha_{jk}\rho_{jk}^4)
 \end{equation}

Taking $\theta=\pi/2$ and with the notation introduced in terms 2, we find the relations: 

$$\rho_{12}=\frac{1}{\sqrt 2 \cos\phi},~\rho_{34}=\frac{1}{\sqrt 2\sin\phi}$$

$$\rho_{13}=\frac{\sqrt 2}{\cos\phi-\sin\phi}=-\rho_{24}$$

$$\rho_{14}=\frac{\sqrt 2}{\cos\phi+\sin\phi}=-\rho_{23}$$

$$\alpha_{12}=\frac{1}{\rho_{34}^2},~\alpha_{34}=\frac{1}{\rho_{12}^2}$$ $$\alpha_{13}=\frac{1}{\rho_{14}^2}+1=\alpha_{24}$$ $$\alpha_{14}=\frac{1}{\rho_{13}^2}+1=\alpha_{23}.$$

Now the left side of (eq. \ref{ineq}) works out to:

$$2((\rho_{12}^3+\rho_{13}^3+\rho_{14}^3)^2+(\rho_{13}^3-\rho_{14}^3-\rho_{34}^3)^2)=$$

$$=2(\sum_{k>j}\rho_{jk}^6+2\rho_{12}^3(\rho_{13}^3+\rho_{14}^3)+2\rho_{34}^3 (\rho_{14}^3-\rho_{13}^3))=2\sum_{k>j}\rho_{jk}^6-96\frac{1}{\sin^2 2\phi\cos^2 2\phi}=$$
$$=2\sum_{k>j}\rho_{jk}^6+\text{negative term}$$

and the right side of (eq. \ref{ineq}) works out to:

$$(\rho_{12}^2+\rho_{34}^2+2(\rho_{13}^2+\rho_{14}^2))(\frac{\rho_{12}^4}{\rho_{34}^2}+\frac{\rho_{34}^4}{\rho_{12}^2}+2(\rho_{13}^4+\rho_{14}^4+\frac{\rho_{13}^4}{\rho_{14}^2}+\frac{\rho_{14}^4}{\rho_{13}^2}))=$$

$$=(\frac{2}{\sin^2 2\phi}+\frac{8}{\cos^2 2\phi})(\sin^2 2\phi (\rho_{12}^6+\rho_{34}^6)+\frac{\cos^2 2\phi}{2}(\rho_{13}^6+\rho_{14}^6)+2(\rho_{13}^4+\rho_{14}^4)=$$

$$=2\sum_{k>j} \rho_{jk}^6+\cot^2 2\phi (\rho_{13}^6+\rho_{14}^6)+8\tan^2 2\phi (\rho_{12}^6+\rho_{34}^6)+(\rho_{13}^4+\rho_{14}^4)(\frac{4}{\sin^2 2\phi}+\frac{16}{\cos^2 2\phi})=$$

$$=2\sum_{k>j}\rho_{jk}^6+\text{positive term}.$$

Therefore the inequality (eq. \ref{ineq}) holds! $\qed$\\

\section{Acknowledgment.}We would like to thank Slobodan Simic [San Jose State, San Jose ,CA] for discussions on hyperbolicity,
and  Prof. Jie Qing [UCSC] and Wei Yuan of UCSC   for conversations regarding curvature computations.
We thankfully acknowledge
support   by NSF grant DMS-20030177.

\begin{appendix}

\section{ DERIVATION OF eq. \ref{curv}:\\}

Take a $ds^2$-orthonormal basis $\{ v_a\}$ for $\C^3$ with $v_1, v_2\in\mathcal{H}_p$.

The Kulkarni-Nomizu [K-N] product formula for conformal curvatures (\cite{KN}, pg. 51) reads:
$$\bar R_{abcd}-U_LR_{abcd}=-\{ d s_{JM}^2 \text{\textcircled{$\wedge$}}(\nabla du -du\otimes du+\frac12 \| du\|^2 ds^2)\}_{abcd}$$
where $u:=\frac12 \log U_L$ and the overbars denote curvature with respect to the $d s_{JM}^2$-metric and all other quantities
(no overbars) are with respect to the $ds^2$-metric.  Then $R_{abcd}=0$ since $ds^2$ is the flat Euclidean metric of  $\C^3 = \R^6$.  Taking $cd=ab$ we have:
$$U_L^2\bar K_{ab}=\bar R_{abab}=-U_L(\nabla du_{bb}+\nabla du_{aa}-du_b\otimes du_b-du_a\otimes du_a+\| du\|^2)=$$
$$=-U_L(\del_a^2 u+\del_b^2 u-(\del_a u)^2-(\del_b u)^2+\|\nabla u\|^2).$$

Next O'Neill's formula (\cite{ON}, pg. 213) gives

$$\mathcal{K}(d\pi v_1, d\pi v_2)=\bar K_{12}+\frac34 |[V_1, V_2]^\mathcal{V}|_{d s_{JM}^2}^2$$

where $V_a=\frac{v_a}{\sqrt{U_L(p)}}$ and $X^{\mathcal{V}}$ denotes $d s_{JM}^2$ projection of $X$ onto $\mathcal{V}$.

We then compute:

$$\del_a u=\frac{\del_a U_L}{2U_L}=\frac{\nabla |_{Lp} U\cdot Lv_a}{2U_L(p)}$$ and 
$$\del_a^2 u=\frac{\del_a^2 U_L}{2U_L}-\frac{(\del_a U_L)^2}{2U_L^2}=\frac{\nabla|_{Lp} (\nabla U\cdot Lv_a)\cdot Lv_a}{2U_L(p)}-\frac{(\del_a U_L)^2}{2U_L(p)^2}.$$

Note that $\nabla U\in \{ q\in \C^4 : \sum q_j=0\}$ and $L v_a$ is a $ds^2$ orthonormal basis for this center of mass zero subspace, hence $$\|\nabla U\|^2=\sum (\nabla U\cdot Lv_a)^2=\sum (\del_a U_L)^2=4U_L^2\|\nabla u\|^2.$$

Substitution into the K-N formula gives 

\begin{equation}
 \label{kn}
\bar K_{12}=-\frac{1}{U_L^3}(\frac{U_L}{2}(\del_1^2 U_L+\del_2^2 U_L)-\frac34 (\del_1 U_L^2+\del_2 U_L^2)+\|\nabla U/2\|^2).
 \end{equation}

To compute O'Neill's Lie bracket term we write our standard coordinates on $\C^3$ as  $(x^1+ ix^2, ..., x^5+ ix^6)$.

Let $H_1=X^j\del_{x^j},~H_2=Y^j\del_{x^j}\in \mathcal{H}$ be any horizontal vector fields.  The vertical vector fields
are spanned by the Euler vector field  $E=x^j\del_{x^j}$ and $iE$. Then $H_j\cdot E=H_j\cdot iE=0$ and:
$$[H_1, H_2]\cdot E=\sum_k X^j x^k\del_{x^j} Y^k-Y^jx^k\del_{x^j}X^k=$$
$$=\sum_k X^j (\del_{x^j} (x^kY^k)-\delta_j^k Y^k) -Y^j(\del_{x^j}(x^kX^k)-\delta_j^k X^k)=\sum_k X^kY^k -Y^k X^k=0$$
and likewise: $$[H_1, H_2]\cdot iE=\sum_{k\text{ odd}} (Y^j\del_{x^j} X^k-X^j\del_{x^j} Y^k)x^{k+1}+(X^j\del_{x^j}Y^{k+1}-Y^j\del_{x^j}X^{k+1})x^k=$$
$$=2\sum_{k\text{ odd}}-X^kY^{k+1}+X^{k+1}Y^k=2H_1\cdot iH_2.$$

Then $$|[V_1,V_2]^{\mathcal{V}_p}|^2=d s_{JM}^2([V_1,V_2], \frac{E_p}{|p|\sqrt{U_L(p)}})^2+d s_{JM}^2([V_1,V_2], \frac{iE_p}{|p|\sqrt{U_L(p)}})^2=$$ 
$$=\frac{U_L^2}{|p|^2U_L}(([V_1,V_2]\cdot E)^2+([V_1,V_2]\cdot iE)^2)=\frac{4U_L(p)(V_1\cdot iV_2)^2}{|p|^2}=\frac{4}{U_L(p)|p|^2}(v_1\cdot iv_2)^2.$$

Now substitution of this Lie bracket expression and (eq. \ref{kn}) into O'Neill's formula and multiplying by $U_L^3$ yields (eq. \ref{curv}). $\qed$

\end{appendix}

 \end{document}